\documentclass[11pt,letterpaper]{article}
\usepackage{charter}
\usepackage{float,amsmath,amsthm,amsfonts,graphicx,color,marvosym}
\usepackage[sort&compress,numbers]{natbib}
\usepackage[lmargin=31mm,rmargin=31mm,tmargin=28mm,bmargin=28mm,footskip=10mm]{geometry}

\newcommand{\seclabel}[1]{\label{sec:#1}}

\newcommand{\secref}[1]{\mbox{Section~\ref{sec:#1}}}
\newcommand{\figlabel}[1]{\label{fig:#1}}

\newcommand{\figref}[1]{\mbox{Figure~\ref{fig:#1}}}

\renewcommand{\eqref}[1]{(\ref{eq:#1})}

\newtheorem{thm}{Theorem}{\bfseries}{\itshape}
\newcommand{\thmlabel}[1]{\label{thm:#1}}
\newcommand{\thmref}[1]{Theorem~\ref{thm:#1}}

\newcommand{\threethmref}[3]{Theorems~\ref{thm:#1}, \ref{thm:#2} and \ref{thm:#3}}

\newtheorem{lem}{Lemma}{\bfseries}{\itshape}
\newcommand{\lemlabel}[1]{\label{lem:#1}}
\newcommand{\lemref}[1]{Lemma~\ref{lem:#1}}
\newcommand{\twolemref}[2]{Lemmas~\ref{lem:#1} and \ref{lem:#2}}

{\bfseries}{\itshape}

\newtheorem{cor}{Corollary}{\bfseries}{\itshape}

\newcommand{\floor}[1]{\left\lfloor #1 \right\rfloor}
\newcommand{\ceil}[1]{\left\lceil #1 \right\rceil}
\newcommand{\Ne}{\ensuremath{\protect{\mathbb{N}}}}
\newcommand{\etal}{et~al.}

\DeclareMathOperator{\dist}{dist}
\newcommand{\e}{\ensuremath{\textup{\textsf{e}}}}
\newcommand{\x}{\ensuremath{\textup{\textsf{x}}}}
\newcommand{\y}{\ensuremath{\textup{\textsf{y}}}}
\newcommand{\dn}[2][]{\ensuremath{\textup{\textsf{dn}}_{#1}(#2)}}
\newcommand{\ddn}[2][]{\ensuremath{\textup{\textsf{ddn}}_{#1}(#2)}}
\newcommand{\DIM}{\ensuremath{\textup{\textsf{dim}}}}
\newcommand{\RAY}[1]{\ensuremath{#1\!\!\downarrow}}
\newcommand{\sn}[2][]{\ensuremath{\textup{\textsf{sn}}_{#1}(#2)}}
\newcommand{\bw}[1]{\textup{\ensuremath{\textsf{bw}(#1)}}}
\newcommand{\cbw}[1]{\textup{\ensuremath{\textsf{cbw}(#1)}}}
\newcommand{\ex}{\textup{\ensuremath{\textsf{ex}}}}
\newcommand{\GG}[1]{\ensuremath{\mathcal{G}\langle{#1}\rangle}}

\title{\MakeUppercase{Distinct Distances in Graph Drawings}}

\author{Paz Carmi\footnotemark[1]
\and Vida Dujmovi\'{c}\footnotemark[1] 
\and Pat Morin\footnotemark[1] 
\and David R.~Wood\footnotemark[2]}

\date{\today}

\renewcommand{\thefootnote}{\arabic{footnote}}

\begin{document}

\maketitle

\renewcommand{\thefootnote}{\fnsymbol{footnote}}

\footnotetext[1]{School of Computer Science, Carleton University, Ottawa, Canada (\texttt{\{paz,vida,morin\}@scs.carleton.ca}). Research supported by NSERC.}

\footnotetext[2]{Departament de Matem{\`a}tica Aplicada II, Universitat Polit{\`e}cnica de Catalunya, Barcelona, Spain (\texttt{david.wood@upc.es}). Research supported by a Marie Curie Fellowship from the European Commission under contract MEIF-CT-2006-023865, and by the projects MEC MTM2006-01267 and DURSI 2005SGR00692.}

\renewcommand{\thefootnote}{\arabic{footnote}}

\begin{abstract}
The \emph{distance-number} of a graph $G$ is the minimum number of distinct edge-lengths over all straight-line drawings of $G$ in the plane. This definition generalises many well-known concepts in combinatorial geometry. We consider the distance-number of trees, graphs with no $K^-_4$-minor, complete bipartite graphs, complete graphs, and cartesian products. Our main results concern the distance-number of graphs with bounded degree. We prove that $n$-vertex graphs with bounded maximum degree and bounded treewidth have distance-number in $\mathcal{O}(\log n)$. To conclude such a logarithmic upper bound, both the degree and the treewidth need to be bounded. In particular, we construct graphs with treewidth $2$ and polynomial distance-number. Similarly, we prove that there exist graphs with maximum degree $5$ and arbitrarily large distance-number. Moreover, as $\Delta$ increases the existential lower bound on the distance-number of $\Delta$-regular graphs tends to $\Omega(n^{0.864138})$.
\end{abstract}

%\newpage
%\tableofcontents
%\newpage

%%%%%%%%%%%%%%%%%%%%%%%%%%%%%%%%%%%%%%%%%%%%%%%%%%
\section{Introduction}\seclabel{intro}
%%%%%%%%%%%%%%%%%%%%%%%%%%%%%%%%%%%%%%%%%%%%%%%%%%

This paper initiates the study of the minimum number of distinct edge-lengths in a drawing of a given graph\footnote{We consider graphs that are simple, finite, and undirected.  The vertex set of a graph $G$ is denoted by $V(G)$, and its edge set by $E(G)$. A graph with $n$ vertices, $m$ edges and maximum degree at most $\Delta$ is an $n$-vertex, $m$-edge, degree-$\Delta$ graph. A graph in which every vertex has degree $\Delta$ is \emph{$\Delta$-regular}. For $S\subseteq V(G)$, let $G[S]$ be the subgraph of $G$ induced by $S$, and let $G-S:=G[V(G)\setminus S]$. For each vertex $v\in V(G)$, let $G-v:=G-\{v\}$. Standard notation is used for graphs: complete graphs $K_n$, complete bipartite graphs $K_{m,n}$, paths $P_n$, and cycles $C_n$. A graph $H$ is a \emph{minor} of a graph $G$ if $H$ can be obtained from a subgraph of $G$ by contracting edges. Throughout the paper, $c$ is a positive constant. Of course, different occurrences of $c$ might denote different constants.}. A \emph{degenerate drawing} of a graph $G$ is a function that maps the vertices of $G$ to distinct points in the plane, and maps each edge $vw$ of $G$ to the open straight-line segment joining the two points representing $v$ and $w$. A \emph{drawing} of $G$ is a degenerate drawing of $G$ in which the image of every edge of $G$ is disjoint from the image of every vertex of $G$. That is, no vertex intersects the interior of an edge. In what follows, we often make no distinction between a vertex or edge in a graph and its image in a drawing. 

%\Comment{Perhaps it would be better to replace ``degenerate drawing'' by ``drawing'', replace ``drawing'' by ``proper drawing'', replace ``degenerate distance-number'' by ``distance-number'', and replace ``distance-number'' by ``proper distance-number''. What do you think?} 

The \emph{distance-number} of a graph $G$, denoted by $\dn{G}$, is the minimum number of distinct edge-lengths in a drawing of  $G$.  The \emph{degenerate distance-number} of $G$, denoted by $\ddn{G}$, is the minimum number of distinct edge-lengths in a degenerate drawing of  $G$. Clearly, $\ddn{G}\leq \dn{G}$ for every graph $G$. Furthermore, if $H$ is a subgraph of $G$ then $\ddn{H}\leq \ddn{G}$ and $\dn{H}\leq \dn{G}$.

%For $v,w\in V(G)$, $\d(v,w)$ denotes the distance between $v$ and $w$ in a (degenerate) drawing of $G$.

%%%%%%%%%%%%%%%%%%%%%%%%%%%%%%%%%%%%%%%%%%%%%%%%%%
\subsection{Background and Motivation}\seclabel{Motivation}
%%%%%%%%%%%%%%%%%%%%%%%%%%%%%%%%%%%%%%%%%%%%%%%%%%

The degenerate distance-number and distance-number of a graph generalise various concepts in combinatorial geometry, which motivates their study. 

A famous problem raised by \citet{Erdos46} asks for the minimum number of distinct distances determined by $n$ points in the plane\footnote{For a detailed exposition on distinct distances in point sets refer to Chapters $10$--$13$ of the monograph by \citet{PA95}.}. This problem is equivalent to determining the degenerate distance-number of the complete graph $K_n$. We have the following bounds on \ddn{K_n}, where the lower bound is due to \citet{KT04} (building on recent advances by \citet{ST01}, \citet{STT-DCG02}, and \citet{Tardos-AM03}), and the upper bound is due to \citet{Erdos46}. 

\begin{lem}[\citep{KT04,Erdos46}]\lemlabel{DegenComplete} 
The degenerate distance-number of $K_n$ satisfies
$$\Omega(n^{0.864137})\leq \ddn{K_n}\leq \frac{cn}{\sqrt{\log n}}\enspace.$$
\end{lem}

Observe that no three points are collinear in a (non-degenerate) drawing of $K_n$. Thus $\dn{K_n}$ equals the minimum number of distinct distances determined by $n$ points in the plane with no three points collinear. This problem was considered by Szemer\'edi (see Theorem 13.7 in \cite{PA95}), who proved that every such point set contains a point from which there are at least $\ceil{\frac{n-1}{3}}$ distinct distances to the other points. Thus we have the next result, where the upper bound follows from the drawing of $K_n$ whose vertices are the points of a regular $n$-gon, as illustrated in \figref{kn}(a). 

\begin{figure}[!ht]
\begin{center}
\includegraphics{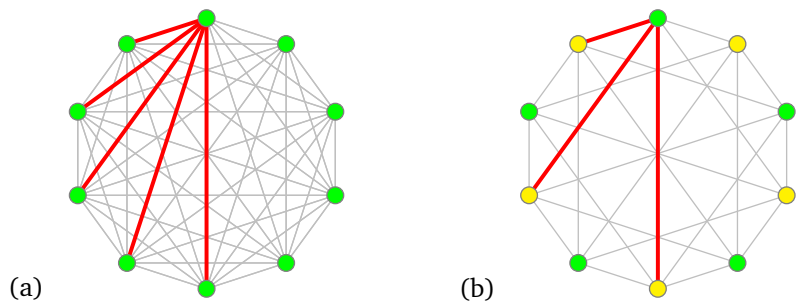}
\caption{\figlabel{kn} (a) A drawing of $K_{10}$ with five edge-lengths, and (b) a drawing of $K_{5,5}$ with three edge-lengths.}\end{center}
\end{figure}

\begin{lem}[Szemer\'edi]\lemlabel{complete}  
The  distance-number of $K_n$ satisfies
$$\ceil{\frac{n-1}{3}}\leq \dn{K_n}\leq \floor{\frac{n}{2}}\enspace.$$
\end{lem}

Note that \twolemref{DegenComplete}{complete} show that for every sufficiently large complete graph, the degenerate distance-number is strictly less than the distance-number. Indeed, $\ddn{K_n}\in o(\dn{K_n})$.

Degenerate distance-number generalises another concept in combinatorial geometry. The \emph{unit-distance graph} of a set $S$ of points in the plane has vertex set $S$, where two vertices are adjacent if and only if they are at unit-distance; see \citep{SheSoi-JCTA03,ODonnell-Geom00b,ODonnell-Geom00,Pritikin-JCTB98,Reid-GC96,HO-Geom96} for example. The famous Hadwiger-Nelson problem asks for the maximum chromatic number of a unit-distance graph. Every unit-distance graph $G$ has $\ddn{G}=1$. But the converse is not true, since a degenerate drawing allows non-adjacent vertices to be at unit-distance. \figref{UnitExample} gives an example of a graph $G$ with $\dn{G}=\ddn{G}=1$ that is not a unit-distance graph. In general, $\ddn{G}=1$ if and only if $G$ is isomorphic to a subgraph of a unit-distance graph. 

\begin{figure}[ht]
\begin{center}
\includegraphics{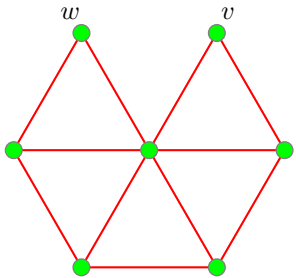}
\caption{\figlabel{UnitExample} A graph with distance-number $1$ that is not a unit-distance graph. In every mapping of the vertices to distinct points in the plane with unit-length edges, $v$ and $w$ are at unit-distance.}\end{center}
\end{figure}

The maximum number of edges in a unit-distance graph is an old open problem. The best construction, due to \citet{Erdos46}, gives an $n$-vertex unit-distance graph with $n^{1+c/\log\log n}$ edges. The best upper bound on the number of edges is $cn^{4/3}$, due to \citet{SST84}. (\citet{Szekely-CPC97} found a simple proof for this upper bound based on the crossing lemma.)\ 

More generally, many recent results in the combinatorial geometry literature provide upper bounds on the number of times the $d$ most frequent inter-point distances can occur between a set of $n$ points. Such results are equivalent to upper bounds on the number of edges in an $n$-vertex graph with degenerate distance number $d$. This suggests the following extremal function. Let $\ex(n,d)$ be the maximum number of edges in an $n$-vertex graph $G$ with $\ddn{G}\leq d$. 

Since every graph $G$ is the union of \ddn{G} subgraphs of unit-distance graphs, the above result by \citet{SST84} implies:

%Take the same grid that gives $n^{1+c/\log\log n}$ bound for $d=1$. 

\begin{lem}[\citet{SST84}]\lemlabel{numedges}
$$\ex(n,d)\leq cdn^{4/3}.$$
Equivalently, the distance-numbers of every $n$-vertex $m$-edge graph $G$ satisfy $$\dn{G}\geq\ddn{G}\geq  cmn^{-4/3}\enspace.$$
\end{lem}

Results by \citet{KT04} (building on recent advances by \citet{ST01}, \citet{STT-DCG02}, and \citet{Tardos-AM03}) imply:

\begin{lem}[\citet{KT04}]\lemlabel{extremal}
$$\ex(n,d)\in\mathcal{O}\big(n^{1.457341}d^{0.627977}\big)\enspace.$$
Equivalently, the distance-numbers of every $n$-vertex $m$-edge graph $G$ satisfy
$$\dn{G}\geq\ddn{G}\in\Omega(m^{1.592412}\,n^{-2.320687})\enspace.$$
\end{lem}

Note that \lemref{extremal} improves upon \lemref{numedges} whenever $\ddn{G}>n^{1/3}$. Also note that \lemref{extremal} implies the lower bound in \lemref{complete}.

\subsection{Our Results}

The above results give properties of various graphs defined with respect to the inter-point distances of a set of points in the plane. This paper, which is more about graph drawing than combinatorial geometry, reverses this approach, and asks for a drawing of a given graph with few inter-point distances. 

Our first results provide some general families of graphs, namely trees and graphs with no $K^-_4$-minor, that are unit-distance graphs (\secref{UnitDistanceGraphs}). Here $K^-_4$ is the graph obtained from $K_4$ by deleting one edge. Then we give bounds on the distance-numbers of complete bipartite graphs (\secref{CompleteBipartite}). 

Our main results concern graphs of bounded degree (\secref{Degree}). We prove that for all $\Delta\geq5$ there are degree-$\Delta$ graphs with unbounded distance-number. Moreover, for $\Delta\geq7$ we prove a polynomial lower bound on the distance-number (of some degree-$\Delta$ graph) that tends to $\Omega(n^{0.864138})$ for large $\Delta$. On the other hand, we prove that graphs with bounded degree and bounded treewidth have distance-number in $\mathcal{O}(\log n)$. Note that bounded treewidth alone does not imply a logarithmic bound on distance-number since $K_{2,n}$ has treewidth $2$ and degenerate distance-number $\Theta(\sqrt{n})$ (see \secref{CompleteBipartite}). 

Then we establish an upper bound on the distance-number in terms of the bandwidth (\secref{Bandwidth}). Then we consider the distance-number of the cartesian product of graphs (\secref{CartesianProducts}). We conclude in \secref{Conclusion} with a discussion of open problems related to distance-number. 

%In what follows we often derive lower bounds for degenerate distance-number and upper bounds for distance-number; then such bounds apply to both quantities. 

\subsection{Higher-Dimensional Relatives}

Graph invariants related to distances in higher dimensions have also been studied. 
\citet*{EHT65} defined the \emph{dimension} of a graph $G$, denoted by $\DIM(G)$, to be the minimum integer $d$ such that $G$ has a degenerate drawing in $\Re^d$ with straight-line edges of unit-length. They proved that $\DIM(K_n)=n-1$, the dimension of the $n$-cube is $2$ for $n \geq 2$, the dimension of the Peterson graph is $2$, and $\DIM(G)\leq 2\cdot\chi(G)$ for every graph $G$. (Here $\chi(G)$ is the \emph{chromatic number} of $G$.)\ The dimension of complete $3$-partite graphs and wheels were determined by \citet{BH-GC88}.

The \emph{unit-distance graph} of a set $S\subseteq\Re^d$ has vertex set $S$, where two vertices are adjacent if and only if they are at unit-distance. Thus $\DIM(G)\leq d$ if and only if $G$ is isomorphic to a subgraph of a unit-distance graph in $\Re^d$. \citet{Mae-DCG89} proved for all $d$ there is a finite bipartite graph (which thus has dimension at most $4$) that is not a unit-distance graph in $\Re^d$. This highlights the distinction between dimension and unit-distance graphs. \citet{Mae-DCG89} also proved that every finite graph with maximum degree $\Delta$ is a unit-distance graph in $\Re^{\Delta(\Delta^2-1)/2}$, which was improved to $\Re^{2\Delta}$ by \citet{MaeRod-GC90}. These results are in contrast to our result that graphs of bounded degree have arbitrarily large distance-number. 

A graph is $d$-\emph{realizable} if, for every mapping of its vertices to (not-necessarily distinct) points in $\Re^p$ with $p\geq d$, there exists such a mapping in $\Re^d$ that preserves edge-lengths. For example, $K_3$ is $2$-realizable but not $1$-realizable. \citet{BC-DCG07} and \citet{Belk-DCG07} proved that a graph is $2$-realizable if and only if it has treewidth at most $2$. They also characterized the $3$-realizable graphs as those with no $K_5$-minor and no $K_{2,2,2}$-minor.

%Relationship to distance-number: Since each graph has a realization in $\Ee^d$ with all edges as unit-distances $-->$ all partial 2-trees have realization in the plane with unit-distances. So that is in contrast to our result that a partial 2-tree, $K_{2,n}$ has distance-number in $\Omega(\sqrt{n})$

%%%%%%%%%%%%%%%%%%%%%%%%%%%%%%%%%%%%%%%%%%%%%%%%%%
\section{Some Unit-Distance Graphs}\seclabel{UnitDistanceGraphs}
%%%%%%%%%%%%%%%%%%%%%%%%%%%%%%%%%%%%%%%%%%%%%%%%%%

This section shows that certain families of graphs are unit-distance graphs. The proofs are based on the fact that two distinct circles intersect in at most two points. We start with a general lemma. A graph $G$ is obtained by \emph{pasting} subgraphs $G_1$ and $G_2$ on a cut-vertex $v$ of $G$ if $G=G_1\cup G_2$ and $V(G_1)\cap V(G_2)=\{v\}$.

\begin{lem}
\lemlabel{Pasting}
Let $G$ be the graph obtained by pasting subgraphs $G_1$ and $G_2$ on a vertex $v$. Then:\\[1ex]
\hspace*{10mm} (a) if $\ddn{G_1}=\ddn{G_2}=1$ then $\ddn{G}=1$, and \\[1ex]
\hspace*{10mm} (b) if $\dn{G_1}=\dn{G_2}=1$ then $\dn{G}=1$.
\end{lem}

\begin{proof}
We prove part (b). Part (a) is easier. 
Let $D_i$ be a drawing of $G_i$ with unit-length edges.
Translate $D_2$ so that $v$ appears in the same position in $D_1$ and $D_2$.
A rotation of $D_2$ about $v$ is \emph{bad} if its union with $D_1$ is not a drawing of $G$. That is, some vertex in $D_2$ coincides with the closure of some edge of $D_1$, or vice versa.
Since $G$ is finite, there are only finitely many bad rotations.
Since there are infinitely many rotations, there exists a rotation that is not bad.
That is, there exists a drawing of $G$ with unit-length edges.
\end{proof}

We have a similar result for unit-distance graphs. 

\begin{lem}
\lemlabel{PastingUnitDistance}
Let $G_1$ and $G_2$ be unit-distance graphs. Let $G$ be the (abstract) graph obtained by pasting $G_1$ and $G_2$ on a vertex $v$. Then $G$ is isomorphic to a unit-distance graph.
\end{lem}

\begin{proof}
The proof is similar to the proof of \lemref{Pasting}, except that we must ensure that the distance between vertices in $G_1-v$ and vertices in $G_2-v$ (which are not adjacent) is not $1$. Again this will happen for only finitely many rotations. Thus there exists a rotation that works.
\end{proof}

Since every tree can be obtained by pasting a smaller tree with $K_2$, \lemref{PastingUnitDistance} implies that every tree is a unit-distance graph. 
The following is a stronger result.

\begin{lem}\lemlabel{trees}
Every tree $T$ has a crossing-free\footnote{A drawing is \emph{crossing-free} if no pair of edges intersect.} drawing in the plane such that two vertices are adjacent if and only if they are unit-distance apart.
\end{lem}

\begin{proof} 
For a point $v=(\x(v),\y(v))$ in the plane, let \RAY{v} be the ray from $v$ to $(\x(v),-\infty)$. We proceed by induction on $n$ with the following hypothesis: Every tree $T$ with $n$ vertices has the desired drawing, such that the vertices have distinct \x-coordinates, and for each vertex $u$, the ray \RAY{u} does not intersect $T$. The statement is trivially true for $n\leq 2$. For $n>2$, let $v$ be a leaf of $T$ with parent $p$. By induction, $T-v$ has the desired drawing. Let $w$ be a vertex of $T-v$, such that no vertex has its \x-coordinate between $\x(p)$ and $\x(w)$. Thus the drawing of $T-v$ does not intersect the open region $R$ of the plane bounded by the two rays \RAY{p} and \RAY{w}, and the segment $pw$. Let $A$ be the intersection of $R$ with the unit-circle centred at $p$. Thus $A$ is a circular arc. Place $v$ on $A$, so that the distance from $v$ to every vertex except $p$ is not $1$. This is possible since $A$ is infinite, and there are only finitely many excluded positions on $A$ (since $A$ intersects a unit-circle centred at a vertex except $p$ in at most two points). Since there are no elements of $T-v$ in $R$, there are no crossings in the resulting drawing and the induction invariants are maintained for all vertices of $T$.
\end{proof}

\begin{figure}[t]
\begin{center}
\includegraphics{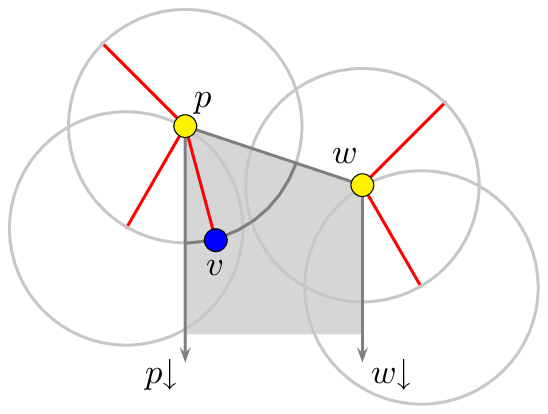}
\end{center}
\caption{Illustration for the proof of \lemref{trees}}
\end{figure}

Recall that $K^-_4$ is the graph obtained from $K_4$ by deleting one edge. 

%Note that a graph is $K^-_4$-minor-free if and only if between each pair of vertices there are at most two internally disjoint paths.

\begin{thm} 
Every $K^-_4$-minor-free graph $G$ has a drawing such that vertices are adjacent if and only if they are unit-distance apart. In particular, $G$ is isomorphic to a unit-distance graph and $\ddn{G}=\dn{G}=1$.
\end{thm}

\begin{proof}
Clearly we can assume that $G$ is connected. If $G$ is not $2$-connected then $G$ can be obtained by pasting smaller graphs with no $K^-_4$-minor on a single vertex. Thus by \lemref{PastingUnitDistance}, it suffices to prove that every $2$-connected $K^-_4$-minor-free graph is a unit-distance graph. 
Below we prove that every $2$-connected graph $G$ with no $K^-_4$-minor is a cycle. The result follows since $C_n$ is a unit-distance graph (draw a regular $n$-gon).

Suppose on the contrary that $G$ has a vertex $v$ of degree at least $3$.
Let $x,y,z$ be the neighbours of $v$.
So there is an $xy$-path $P$ avoiding $v$ (since $G$ is $2$-connected) and avoiding $z$ (since $G$ is $K^-_4$-minor free).
Similarly, there is an $xz$-path $Q$ avoiding $v$.
If $x$ is the only vertex in both $P$ and $Q$, 
then the cycle $(x,P,y,v,z,Q)$ plus the edge $xv$ is a subdivision of $K^-_4$.
Now assume that $P$ and $Q$ intersect at some other vertex.
Let $t$ be the first vertex on $P$ starting at $x$ that is also in $Q$.
Then the cycle $(x,Q,z,v)$ plus the sub-path of $P$ between $x$ and $t$ is a subdivision of $K^-_4$. 
Hence $G$ has no vertex of degree at least $3$. 
Therefore $G$ is a cycle, as desired.
\end{proof}

%\Comment{Suppose that for some graph $H$, every $H$-minor-free graph is a unit-distance graph. Conjecture: $H\subseteq K^-_4$.}

%%%%%%%%%%%%%%%%%%%%%%%%%%%%%%%%%%%%%%%%%%%%%%%%%%%%%%%%%%%%%%%%%%%
\section{Complete Bipartite Graphs}\seclabel{CompleteBipartite}
%%%%%%%%%%%%%%%%%%%%%%%%%%%%%%%%%%%%%%%%%%%%%%%%%%%%%%%%%%%%%%%%%%%

This section considers the distance-numbers of complete bipartite graphs $K_{m,n}$. Since $K_{1,n}$ is a tree, $\ddn{K_{1,n}}=\dn{K_{1,n}}=1$ by \lemref{trees}. The next case, $K_{2,n}$, is also easily handled.

\begin{lem}\lemlabel{k2n}
The distance-numbers of $K_{2,n}$ satisfy
$$\ddn{K_{2,n}}=\dn{K_{2,n}} = \ceil{\sqrt\frac{n}{2}}\enspace.$$
\end{lem}

\begin{proof}
Let $G=K_{2,n}$ with colour classes $A=\{v,w\}$ and $B$, where $|B|=n$. We first prove the lower bound. Consider a degenerate drawing of $G$ with $\ddn{G}$ edge-lengths. The vertices in $B$ lie on the intersection of $\ddn{G}$ concentric circles centered at $v$ and $\ddn{G}$ concentric circles centered at $w$. 
Since two distinct circles intersect in at most two points, $n\leq 2\,\ddn{G}^2$. Thus $\ddn{K_{2,n}}\geq\ceil{\sqrt\frac{n}{2}}$

% Thus one of these circles, $C$, contains at least $\frac{n}{\ddn{G}}$ vertices. Since two circles intersect in at most two points, the edges between $w$ and the vertices in $V\cap C$ determine at least $\frac{n}{2\ddn{G}}$ distances. Thus $\ddn{G}\geq \frac{n}{2\ddn{G}}$, giving $\ddn{G}\geq \ceil{\sqrt\frac{n}{2}}$.

For the upper bound, position $v$ at $(-1,0)$ and $w$ at $(1,0)$. As illustrated in \figref{k2n}, draw $\ceil{\sqrt\frac{n}{2}}$ circles centered at each of $v$ and $w$ with radii ranging strictly between $1$ and $2$, such that the intersections of the circles together with $v$ and $w$ define a set of points with no three points collinear. (This can be achieved by choosing the radii iteratively, since for each circle $C$, there are finitely many forbidden values for the radius of $C$.) Each pair of non-concentric circles intersect in two points. Thus the number of intersection points is at least $n$. Placing the vertices of $B$ at these intersection points results in a drawing of $K_{2,n}$ with $\ceil{\sqrt\frac{n}{2}}$ edge-lengths.\end{proof}

\begin{figure}[ht]
\begin{center}
\includegraphics{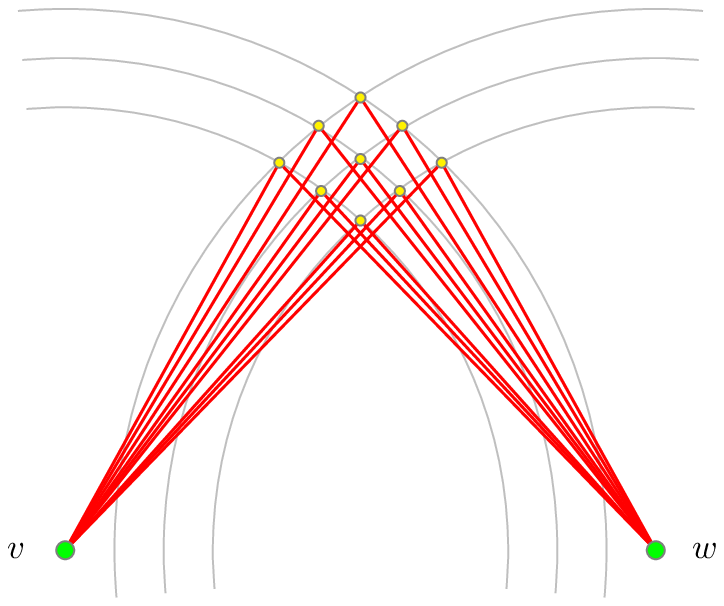}
\caption{\figlabel{k2n} Illustration for the proof of \lemref{k2n}.}\end{center}
\end{figure}

Now we determine \ddn{K_{3,n}} to within a constant factor.

\begin{lem}\lemlabel{k3n}
The degenerate distance-number of $K_{3,n}$ satisfies
$$\ceil{\sqrt\frac{n}{2}}\leq \ddn{K_{3,n}}\leq  3\ceil{\sqrt\frac{n}{2}}-1\enspace.$$
\end{lem}

%\Comment{Is $\dn{K_{3,n}}\leq c\sqrt{n}$? The same construction with no vertices in $B$ at $\x=0$ might work???}

\begin{proof}
The lower bound follows from \lemref{k2n} since $K_{2,n}$ is a subgraph of $K_{3,n}$.

Now we prove the upper bound. Let $A$ and $B$ be the colour classes of $K_{3,n}$, where $|A|=3$ and $|B|=n$. Place the vertices in $A$ at $(-1,0)$, $(0,0)$, and $(1,0)$. Let $d:=\ceil{\sqrt\frac{n}{2}}$. For $i\in[d]$, let $$r_i:=\sqrt{1+\frac{i}{d+1}}\enspace.$$ Note that $1<r_i<2$. 
Let $R_i$ be the circle centred at $(-1,0)$ with radius $r_i$. 
For $j\in[d]$, let $S_j$ be the circle centred at $(1,0)$ with radius $r_j$. 
Observe that each pair of circles $R_i$ and $S_j$ intersect in exactly two points. Place the vertices in $B$ at the intersection points of these circles. 
This is possible since $2d^2\geq n$.

Let $(x,y)$ and $(x,-y)$ be the two points where $R_i$ and $S_j$ intersect. 
Thus $(x+1)^2+y^2=r_i^2$ and $(x-1)^2+y^2=r_j^2$. 
It follows that $$x^2+y^2=\frac{i}{d+1}+2x=\frac{j}{d+1}-2x\enspace.$$ 
Thus $2(x^2+y^2)=\frac{i+j}{d+1}$. That is, the distance from $(x,y)$ to $(0,0)$ equals $$\sqrt{\frac{i+j}{2d+2}},$$ which is the same distance from $(x,-y)$ to $(0,0)$.
Thus the distance from each vertex in $B$ to $(0,0)$ is one of $2d-1$ values  (determined by $i+j$). The distance from each vertex in $B$ to $(-1,0)$ and to $(1,0)$ is one of $d$ values. Hence the degenerate distance-number of $K_{3,n}$ is at most $3d-1=3\ceil{\sqrt\frac{n}{2}}-1$.
\end{proof}

\begin{figure}[ht]
\begin{center}
\includegraphics{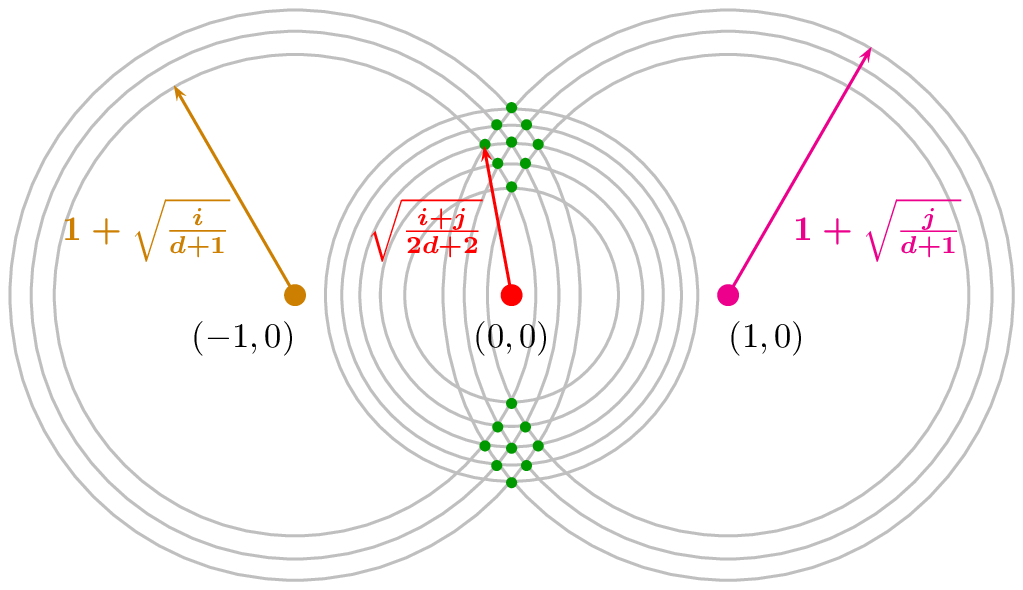}
\caption{\figlabel{k3n} Illustration for the proof of \lemref{k3n}.}\end{center}
\end{figure}

Now consider the distance-number of a general complete bipartite graph.

\begin{lem}\lemlabel{kmn}
For all $n\geq m$, the distance-numbers of $K_{m,n}$ satisfy
$$\Omega\Big(\frac{mn}{(m+n)^{1.457341}}\Big)^{1/0.627977}\leq 
\ddn{K_{m,n}} \leq \dn{K_{m,n}} \leq \ceil{\frac{n}{2}}\enspace.$$
In particular, 
$$\Omega(n^{0.864137})\leq \ddn{K_{n,n}} \leq \dn{K_{n,n}} \leq \ceil{\frac{n}{2}}\enspace.$$
\end{lem}

\begin{proof}
The lower bounds follow from \lemref{extremal}. For the upper bound on $\dn{K_{n,n}}$, position the vertices on a regular $2n$-gon $(v_1,v_2,\dots,v_{2n})$ alternating between the colour classes, as illustrated in \figref{kn}(b). In the resulting drawing of $K_{n,n}$, the number of edge-lengths is $|\{(i+j)\mod n\,:\, v_iv_j\in E(K_{n,n})\}|$. Since $v_iv_j$ is an edge if and only if $i+j$ is odd, the number of edge-lengths is $\ceil{\frac{n}{2}}$. The upper bound on \dn{K_{n,m}} follows since $K_{n,m}$ is a subgraph of $K_{n,n}$.
\end{proof}

%%%%%%%%%%%%%%%%%%%%%%%%%%%%%%%%%%%%%%%%%%%%%%%%%%%%%%%%%%%%%%%%
\section{Bounded degree graphs}\seclabel{Degree}
%%%%%%%%%%%%%%%%%%%%%%%%%%%%%%%%%%%%%%%%%%%%%%%%%%%%%%%%%%%%%%%%

\lemref{k2n} implies that if a graph has two vertices with many common neighbours then its distance-number is necessarily large. Thus it is natural to ask whether graphs of bounded degree have bounded distance-number. This section provides a negative answer to this question. 

%%%%%%%%%%%%%%%%%%%%%%%%%%%%%%%%%%%%%%%%%%%%%%%%%%
\subsection{\boldmath Bounded degree graphs with $\Delta\geq 7$}\seclabel{7}
%%%%%%%%%%%%%%%%%%%%%%%%%%%%%%%%%%%%%%%%%%%%%%%%%%

%The following theorem is the main result of this section.

%%%%%%%%%%%%%%%%%%%%%%%%
%\begin{thm}\thmlabel{degree-7}
%For all $\Delta\geq 7$ and for all $\varepsilon>0$, for all sufficiently large $n>n(\varepsilon, \Delta)$, there exists a $\Delta$-regular $n$-vertex graph $G$ with degenerate distance-number
%$$\ddn{G}>n^{\frac{2}{3}-\frac{20+10\varepsilon}{3\Delta+12+6\varepsilon}}\enspace.$$
%\end{thm}
%%%%%%%%%%%%%%%%%%%%%%%%

This section proves that for all $\Delta\geq7$ there are $\Delta$-regular graphs with unbounded distance-number. Moreover, the lower bound on the distance-number is polynomial in the number of vertices. The basic idea of the proof is to show that there are more $\Delta$-regular graphs than graphs with bounded distance-number. 

It will be convenient to count labelled graphs. Let \GG{n,\Delta} denote the family of labelled $\Delta$-regular $n$-vertex graphs. Let \GG{n,m,d} denote the family of labelled $n$-vertex $m$-edge graphs with degenerate distance-number at most $d$. Our results follow by comparing a lower bound on $|\GG{n,\Delta}|$ with an upper bound on $|\GG{n,m,d}|$ with $m=\frac{\Delta n}{2}$, which is the number of edges in a $\Delta$-regular $n$-vertex graph.

The lower bound in question is known. In particular, the first asymptotic bounds on the number of labelled $\Delta$-regular $n$-vertex graphs were independently determined by \citet{BC-JCTA78} and \citet{Wormald78}. \citet{McKay-AC85} further refined these results. We will use the following simple lower bound derived by \citet{BMW-EJC06} from the result of \citet{McKay-AC85}.

%%%%%%%%%%%%%%%%%%%%%%%%
\begin{lem}[\citep{BMW-EJC06,BC-JCTA78,Wormald78,McKay-AC85}]
\lemlabel{NumberRegular}
For all integers $\Delta\geq1$ and $n\geq c\Delta$, the number of labelled $\Delta$-regular $n$-vertex graphs satisfies
\begin{equation*}
|\GG{n,\Delta}|\geq\Big(\frac{n}{3\Delta}\Big)^{\Delta n/2}.
\end{equation*}
\end{lem} 
%%%%%%%%%%%%%%%%%%%%%%%%

The proof of our upper bound on $|\GG{n,m,d}|$ uses the following special case of the Milnor-Thom theorem by \citet{RBG-JAMS01}. Let $\mathcal{P}=(P_1,P_2,\dots,P_t)$ be a sequence of polynomials on $p$ variables over $\Re$. The \emph{zero-pattern} of ${\mathcal P}$ at $u\in \Re^p$ is the set $\{i:1\leq i\leq t,\,P_i(u)=0\}$. 

%%%%%%%%%%%%%%%%%%%%%%%%
\begin{lem}[\citep{RBG-JAMS01}]\lemlabel{zeros}
Let ${\mathcal P} = (P_1, P_2, \dots, P_t)$ be a sequence of polynomials of degree at most $\delta\geq 1$ on $p\leq t$ variables over $\Re$. Then the number of zero-patterns of ${\mathcal P}$ is at most $\binom{\delta t}{p}$.
\end{lem}

Recall that $\ex(n,d)$ is the maximum number of edges in an $n$-vertex graph $G$ with $\ddn{G}\leq d$. Bounds on this function are given in \twolemref{numedges}{extremal}. Our upper bound on $|\GG{n,m,d}|$ is expressed in terms of $\ex(n,d)$.

%%%%%%%%%%%%%%%%%%%%%%%%
\begin{lem}\lemlabel{num}
The number of labelled $n$-vertex $m$-edge graphs with $\ddn{G}\leq d$ satisfies
$$|\GG{n,m,d}|\leq\left(\frac{\e nd}{2}\right)^{2n+d}\binom{\ex(n,d)}{m}\enspace,$$
where \e\ is the base of the natural logarithm. 
\end{lem}

\begin{proof}
Let $V(G)=\{1,2,\dots,n\}$ for every $G\in\GG{n,m,d}$. For every $G\in\GG{n,m,d}$, there is a point set $$S(G)=\{(x_i(G), y_i(G))\,:\, 1\leq i \leq n\}$$ and a set of edge-lengths $$L(G) = \{\ell_k(G)\,:\, 1\leq k\leq d\},$$ such that $G$ has a degenerate drawing in which each vertex $i$ is represented by the point $(x_i(G), y_i(G))$ and the length of each edge in $E(G)$ is in $L(G)$. Fix one such degenerate drawing of $G$. 

For all $i,j,k$ with $1\leq i<j\leq n$ and $1\leq k\leq d$, and for every graph $G\in\GG{n,m,d}$, define
$$P_{i,j,k}(G) := (x_j(G)-x_i(G))^2+(y_j(G)-y_i(G))^2- \ell_k(G)^2\enspace.$$
Consider ${\mathcal P}:= \{P_{i,j,k}: 1\leq i<j\leq n, 1\leq k\leq d\}$ to be a set of $\binom{n}{2}d$  degree-$2$ polynomials on the set of $2n+d$ variables $\{x_1, x_2,\dots, x_n, y_1, y_2,\dots, y_n, \ell_1, \ell_2,\dots,\ell_d\}$.
Observe that

\smallskip
\begin{minipage}[t]{12cm}  
$P_{i,j,k}(G)=0$ if and only if the distance between vertices $i$ and $j$ in the\\ degenerate drawing of $G$ is $\ell_k(G)$.
\end{minipage}
\hspace*{3mm}
$(\star)$

\smallskip
%if $ij$ is an edge of $G$ and $ij$ has edge-length $\ell_k$ in the representation of $G$ then $P_{i,j,k}(G)=0$. 
Recall the  well-known fact that $\binom{a}{b}\leq (\frac{\e a}{b})^b$. By \lemref{zeros} with $t=\binom{n}{2}d$, $\delta=2$ and $p=2n+d$, the number of zero-patterns determined by ${\mathcal P}$ is at most
$$
\binom{2\binom{n}{2}d}{2n+d}\leq 
\left(\frac{2\e\binom{n}{2} d}{2n+d} \right)^{2n+d}<
\left(\frac{\e n^2 d}{2n+d} \right)^{2n+d}<
\left(\frac{\e n^2 d}{2n} \right)^{2n+d}= 
\left(\frac{\e nd}{2} \right)^{2n+d}\enspace.$$

Fix a zero-pattern $\sigma$ of $\mathcal{P}$. Let $\mathcal{G}_\sigma$ be the set of graphs $G$ in \GG{n,m,d} such that $\sigma$ is the zero-pattern of ${\mathcal P}$ evaluated at $G$. To bound $|\GG{n,m,d}|$ we now bound $|{\mathcal G}_\sigma|$. Let $H_\sigma$ be the graph with vertex set $V(H_\sigma)=\{1, \dots, n\}$ and edge set $E(H_\sigma)$ where $ij\in E(H_\sigma)$ if and only if $ij\in E(G)$ for some $G\in\mathcal{G}_\sigma$. Consider a degenerate drawing of an arbitrary graph $G\in\mathcal{G}_\sigma$ on the point set $S(G)$. By $(\star)$, $S(G)$ and $L(G)$ define a degenerate drawing of $H$ with $d$ edge-lengths. Thus $\ddn{H_\sigma}\leq d$ and by assumption, $|E(H_\sigma)|\leq \ex(n,d)$. Since every graph in $\mathcal{G}_\sigma$ is a subgraph of $H_\sigma$,  $|\mathcal{G}_\sigma|\leq\binom{|E(H_\sigma)|}{m}$. Therefore, 
$$|\GG{n,m,d}|\leq \left(\frac{\e nd}{2}\right)^{2n+d} \binom{|E(H_\sigma)|}{m} 
\leq \left(\frac{\e nd}{2}\right)^{2n+d}\binom{\ex(n,d)}{m},$$  as required.
\end{proof}

By comparing the lower bound in \lemref{NumberRegular} and the upper bound in \lemref{num} we obtain the following result.

%%%%%%%%%%%%%%%%%%%%%%%%%%%%%%%%%%%%%%%%%%%%%%%%%%%%%%%%%%%%%%%%%%%%%%%%%%%%%%%%%
\begin{lem}\lemlabel{PolyBound}
Suppose that for some real numbers $\alpha$ and $\beta$ with $\beta>0$ and
$1<\alpha<2<\alpha+\beta$, $$\ex(n,d)\in\mathcal{O}(n^{\alpha}d^{\beta})\enspace.$$ Then for every integer $\Delta>\frac{4}{2-\alpha}$, for all $\varepsilon>0$, and for all sufficiently large $n>n(\alpha,\beta,\Delta,\varepsilon)$, there exists a $\Delta$-regular $n$-vertex graph $G$ with degenerate distance-number
$$\ddn{G}>n^{\frac{2-\alpha}{\beta}-
\frac{(2-\alpha+\beta)(4+2\varepsilon)}{\beta^2\Delta+4\beta}} \enspace.$$
\end{lem}

\begin{proof}
In this proof, $\alpha$, $\beta$, $\Delta$ and $\epsilon$ are fixed numbers satisfying the assumptions of the lemma. Let $d$ be the maximum degenerate distance number of a graph in $\GG{n,\Delta}$. The result will follow by showing that for all sufficiently large $n>n(\alpha,\beta,\Delta,\varepsilon)$,
$$d>n^{\frac{2-\alpha}{\beta}-
\frac{(2-\alpha+\beta)(4+2\varepsilon)}{\beta^2\Delta+4\beta}}\enspace.$$
By the definition of $d$, and since every $\Delta$-regular $n$-vertex graph has $\frac{\Delta n}{2}$ edges, every graph in 
\GG{n,\Delta} is also in \GG{n,\frac{\Delta n}{2},d}. By \lemref{NumberRegular} with $n\geq c\Delta$, and by \lemref{num}, 
\begin{align*}
\Big(\frac{n}{3\Delta}\Big)^{\Delta n/2} 
\leq |\GG{n,\Delta}|
\leq |\GG{n,\frac{\Delta n}{2},d}|
\leq \left(\frac{\e nd}{2}\right)^{2n+d}\binom{\ex(n,d)}{\Delta n/2}
\enspace.
\end{align*}
Since $\ex(n,d)\in\mathcal{O}(n^\alpha d^\beta)$, and since $d$ is a function of $n$, there is a constant $c$ such that $\ex(n,d)\leq cn^{\alpha}d^{\beta}$ for sufficiently large $n$. Thus (and since $\binom{a}{b}\leq (\frac{\e a}{b})^b$), 
\begin{align*}
\Big(\frac{n}{3\Delta}\Big)^{\Delta n/2} 
\leq \left(\frac{\e nd}{2}\right)^{2n+d}\binom{ cn^{\alpha}d^{\beta} }{\Delta n/2}
\leq\left(\frac{\e nd}{2}\right)^{2n+d} \left(\frac{ 2\e cn^{\alpha}d^{\beta} }{\Delta n}\right)^{\Delta n/2}\enspace.
\end{align*}
Hence
\begin{align*}
n^{\Delta n} 
\leq 3^{\Delta n} 
\left(\frac{\e nd}{2}\right)^{4n+2d} \left(2\e cn^{\alpha-1}d^{\beta} \right)^{\Delta n}\enspace.
\end{align*}
By \lemref{complete}, $d\leq\ddn{K_n}\leq\frac{cn}{\sqrt{\log n}}$, implying $2d\leq\varepsilon n$ for all large $n>n(\varepsilon)$. Thus
\begin{align*}
n^{\Delta} 
&\leq 3^{\Delta} 
\left(\frac{\e nd}{2}\right)^{4+\varepsilon} \left(2\e cn^{\alpha-1}d^{\beta} \right)^{\Delta}\enspace.
\end{align*}
Hence
\begin{align*}
n^{(2-\alpha)\Delta-4-\varepsilon}
&\leq 3^{\Delta} 
\left(\frac{\e}{2}\right)^{4+\varepsilon}
\left(2\e c\right)^{\Delta}
d^{\beta\Delta+4+\varepsilon}
\enspace.
\end{align*}
Observe that $3^{\Delta}\left(\frac{\e}{2}\right)^{4+\varepsilon} (2\e c)^{\Delta} \leq n^{\varepsilon}$ for all large $n>n(\Delta,\varepsilon)$. Thus
\begin{align*}
n^{(2-\alpha)\Delta-4-2\varepsilon} &\leq d^{\beta\Delta+4+\varepsilon}
\end{align*}
Hence
\begin{align*}
d 
\geq n^{\frac{(2-\alpha)\Delta-4-2\varepsilon}{\beta\Delta+4+\varepsilon}} 
= n^{\frac{2-\alpha}{\beta}-
\frac{(2-\alpha+\beta)(4+\varepsilon)+\beta\epsilon}{\beta(\beta\Delta+4+\varepsilon)}} 
> n^{\frac{2-\alpha}{\beta}-
\frac{(2-\alpha+\beta)(4+2\varepsilon)}{\beta^2\Delta+4\beta}}
\enspace,
\end{align*}
as required.
%This contradiction completes the proof.
\end{proof}

We can now state the main results of this section. By \lemref{numedges}, the conditions of \lemref{PolyBound} are satisfied with $\alpha=\frac{4}{3}$ and $\beta=1$; thus together they imply:

\begin{thm}\thmlabel{degree-7}
For every integer $\Delta\geq 7$, for all $\varepsilon>0$, and for all sufficiently large $n>n(\Delta,\varepsilon)$, there exists a $\Delta$-regular $n$-vertex graph $G$ with degenerate distance-number
$$\ddn{G}>n^{\frac{2}{3}-\frac{20+10\varepsilon}{3\Delta+12}} 
\enspace.$$
\end{thm}

By \lemref{extremal}, the conditions of \lemref{PolyBound} are satisfied with
$\alpha=1.457341$ and $\beta=0.627977$; thus together they imply:

\begin{thm}\thmlabel{degree-8}
For every integer $\Delta\geq8$, for all $\varepsilon>0$, and for all sufficiently large $n>n(\Delta,\varepsilon)$, there exists a $\Delta$-regular $n$-vertex graph $G$ with degenerate distance-number
$$\ddn{G}>
n^{0.864138-\frac{4.682544+2.341272\varepsilon}{0.394355\Delta+2.511908}} \enspace.$$
\end{thm}

Note that the bound given in \thmref{degree-8} is better than the bound in \thmref{degree-7} for $\Delta\geq17$.

%%%%%%%%%%%%%%%%%%%%%%%%%%%%%%%%%%%%%%%%%%%%%%%%%%
\subsection{Bounded degree graphs with $\Delta\geq 5$}\seclabel{5}
%%%%%%%%%%%%%%%%%%%%%%%%%%%%%%%%%%%%%%%%%%%%%%%%%%

\thmref{degree-7} shows that for $\Delta\geq7$ and for sufficiently large $n$, there is an $n$-vertex degree-$\Delta$ graph whose degenerate 
distance-number is at least polynomial in $n$. We now prove that the degenerate distance-number of degree-$5$ graphs can also be arbitrarily large. However, the lower bound we obtain in this case is polylogarithmic in $n$. The proof is inspired by an analogous proof about the slope-number of degree-$5$ graphs, due to \citet{PachPal-EJC06}.

\begin{thm}\thmlabel{degree-5}
For all $d\in \Ne$, there is a degree-$5$ graph $G$ with degenerate distance-number $\ddn{G}>d$.
\end{thm}

\begin{proof}
%Throughout this proof $c>0$ is a constant that possibly changes from one equation to the next. 
Consider the following degree-$5$ graph $G$. For $n\equiv 0\pmod{6}$, let $F$ be the graph with vertex set $\{v_1, v_2,\dots, v_n\}$ and edge set $\{v_iv_j\,:\, |i-j|\leq 2\}$. Let $S:=\{v_i\,:\, i \equiv 1 \pmod{3} \}$. No pair of vertices in $S$ are adjacent in $F$, and $|S|=\frac{n}{3}$ is even. 

Let $\mathcal{M}$ denote the set of all perfect matchings on $S$. For each perfect matching $M_k\in \mathcal{M}$, let $G_k:=F\cup M_k$. As illustrated in \figref{Frame}, let $G$ be the disjoint union of all the $G_k$. Thus the number of connected components of $G$ is $|\mathcal{M}|$, which is at least $(\frac{n}{9})^{n/6}$ by \lemref{NumberRegular} with $\Delta=1$. Here we consider perfect matchings to be $1$-regular graphs. (It is remarkable that even with $\Delta=1$, \lemref{NumberRegular} gives such an accurate bound, since the actual number of matchings in $S$ is $\sqrt{2}(\frac{n}{3\e})^{n/6}$ ignoring lower order additive terms\footnote{For even $n$, let $f(n)$ be the number of perfect matchings of $[n]$. Here we determine the asymptotics of $f$. In every such matching, $n$ is matched with some number in $[n-1]$, and the remaining matching is isomorphic to a perfect matching of $[n-2]$. Every matching obtained in this way is distinct. Thus $f(n)=(n-1)\cdot f(n-2)$, where $f(2)=1$. Hence $f(n)=(n-1)!!=(n-1)(n-3)(n-5)\dots 1$, where $!!$ is the double factorial function. Now $(2n-1)!!=\frac{(2n)!}{2^nn!}$. Thus $f(n)=\frac{n!}{2^{n/2}(n/2)!}\approx\sqrt{2}\,(\frac{n}{\e})^{n/2}$ by Stirling's Approximation.}.)\ 

%=\frac{\sqrt{2\pi{}n}\,(n/\e)^n}{2^{n/2}\,\sqrt{\pi n}\,(n/2\e)^{n/2}} 

\begin{figure}[ht]
\begin{center}
\includegraphics{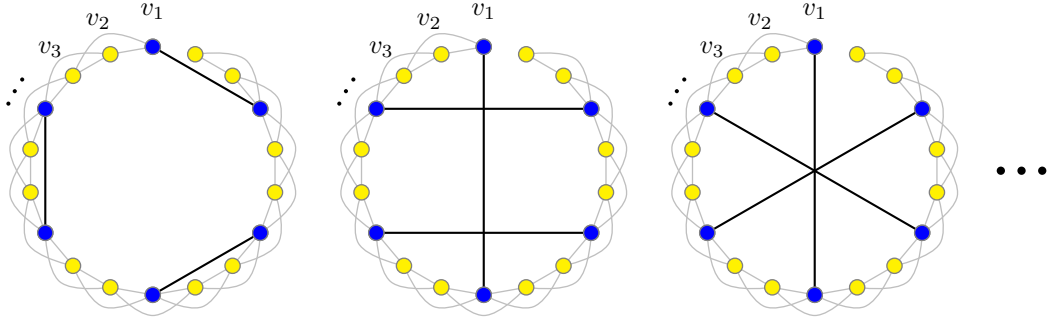}
\caption{\figlabel{Frame} The graph $G$ with $n=18$.}
\end{center}
\end{figure}

Suppose, for the sake of contradiction, that for some constant $d$, for all $n\in \Ne$ such that $n\equiv 0\pmod{6}$, $G$ has a degenerate drawing $D$ with at most $d$ edge-lengths. 

Label the edges of $G$ that are in the copies of $F$ by their length in $D$. Let $\ell_k(i,j)$ be the label of the edge $v_iv_j$ in the copy of $F$ in the component $G_k$ of $G$. This defines a labelling of the components of $G$. Since $F$ has $2n-3$ edges and each edge in $F$ receives one of $d$ labels, there are at most $d^{2n-3}$ distinct labellings of the components of $G$. 

Let $D_k$ be the degenerate drawing of $G_k$ obtained from $D$ by a translation and rotation so that $v_1$ is at $(0,0)$ and $v_2$ is at $(\ell_k(1,2),0)$.
%, and if $i\geq3$ is the minimum integer such that $v_i$ is not on the \x-axis, then $v_i$ has a positive \y-coordinate.
We say that two components $G_q$ and $G_r$ of $G$ \emph{determine the same set of points} if for all $i\in[n]$, the vertex $v_i$ in $D_q$ is at the same position as the vertex $v_i$ in $D_r$. 

Partition the components of $G$ into the minimum number of parts such that all the components in each part have the same labelling and determine the same set of points.

Observe that two components of $G$ with the same labelling do not necessarily determine the same set of points. However, the number of point sets determined by the components with a given labelling can be bounded as follows. For each component $G_k$ of $G$, $v_1$ is at $(0,0)$ and $v_2$ is at $(\ell_k(i,j),0)$ in $D_k$. Thus for a fixed labelling, the positions of $v_1$ and $v_2$ in $D_k$ are determined. Now for $i\geq3$, in each component $G_k$, the vertex $v_i$ is positioned in $D_k$ at the intersection of the circle of radius $\ell_k(i-1,i)$ centered at $v_{i-1}$ and the circle of radius $\ell_k(i-2,i)$ centered at $v_{i-2}$. Thus there are at most two possible locations for $v_i$ (for a fixed labelling). Hence the components with the same labelling determine at most $2^{n-2}$ distinct points sets. Therefore the number of parts in the partition is at most $d^{2n-3}\cdot 2^{n-2}< (2d^2)^n$. 

Finally, we bound the number of components in each part, $R$, of the partition.
Let $H_R$ be the graph with vertex set $V(H_R)=\{v_1,\dots,v_n\}$ where $v_iv_j\in E(H_R)$ if and only if $v_iv_j\in E(G_k)$ for some component $G_k\in R$. Since the graphs in $R$ determine the same set of points, the union of the degenerate drawings $D_k$, over all $G_k\in R$, determines a degenerate drawing of $H_R$ with $d$ edge-lengths. Thus $\ddn{H_R}\leq d$ and by \lemref{numedges}, $|E(H_R)|\leq cdn^{4/3}$ for some constant $c>0$. Every component in $R$ is a subgraph of $H_R$, and any two components in $R$ differ only by the choice of a matching on $S$. Each such matching has $\frac{n}{6}$ edges. Thus the number of components of $G$ in $R$ is at most 
$$\binom{|E(H_R)|}{n/6}
\leq\binom{cdn^{4/3}}{n/6}
\leq\left(\frac{\e cdn^{4/3}}{n/6}\right)^{n/6}
\leq \big(6\e cd\big)^{n/6}n^{n/18}\enspace.$$ 
Hence $|\mathcal{M}|< (2d^2)^n \cdot (6\e cd)^{n/6}n^{n/18}$, and by the lower bound on $|\mathcal{M}|$ from the start of the proof,
$$\left(\frac{n}{9}\right)^{n/6} \leq |\mathcal{M}|
< (2d^2)^n \cdot (6\e cd)^{n/6}n^{n/18}\enspace.$$
The desired contradiction follows for all $n\geq(3456\e c d^{13})^{3/2}$.
\end{proof}

%This assumption implies that the positions of the vertices in $D_k$ are determined by the edge-lengths $\ell_k(i,j)$ (since for $i\geq3$, the vertices $v_{i-2},v_{i-1},v_i$ are a triangle in $G_k$, and the position of the third vertex in a triangle is determined by the position of the first two vertices plus the edge lengths). Note that $D_k$ has at most $d$ edge-lengths. 

%We say that two components $G_q$ and $G_r$ of $G$ \emph{determine the same set of labelled points}, if there exists a point set $P=\{1,\dots,n\}$ such that  for all $1\leq i\leq n$, $v_i\in V(G_q)$ is represented by $i$ in the degenerate drawing of $G_q$, and $v_i\in V(G_r)$ is represented by $i$ in the degenerate drawing of $G_r$.

%; that is, it uniquely determines a set of labelled points.
%We say that two components $G_p$ and $G_q$ of $G$ \emph{determine the same set of labelled points}, if for all $1\leq i\leq n$, $v_i\in V(G_p)$ is represented by the same point in the degenerate drawing of $G_p$ as the vertex $v_i\in V(G_q)$ in the degenerate drawing of $G_q$.
%vertex $v_i$ of $G_p$ and $v_i$ of $G_q$ are represented by the same point in their respective degenerate drawings.

%and  every vertex $v_i$, $1\leq i\leq n$, is represented by the same point in each of the degenerate drawings in that part. Since, $F$ has at most $2n$ edges, there are at most $d^{2n}$ distinct labellings of the components of $G$. 

%%%%%%%%%%%%%%%%%%%%%%%%%%%%%%%%%%%%%%%%%%%%%%%%%%
\subsection{Graphs with bounded degree and bounded treewidth}\seclabel{tw}
%%%%%%%%%%%%%%%%%%%%%%%%%%%%%%%%%%%%%%%%%%%%%%%%%

This section proves a logarithmic upper bound on the distance-number of graphs with bounded degree and bounded treewidth. Treewidth is an important parameter in Robertson and Seymour's theory of graph minors and in algorithmic complexity (see the surveys \cite{Bodlaender-TCS98, Reed-AlgoTreeWidth03}). It can be defined as follows. A graph $G$ is a \emph{$k$-tree} if either $G=K_{k+1}$, or $G$ has a vertex $v$ whose neighbourhood is a clique of order $k$ and $G-v$ is a $k$-tree. For example, every $1$-tree is a tree and every tree is a $1$-tree. Then the \emph{treewidth} of a graph $G$ is the minimum integer $k$ for which $G$ is a subgraph of a $k$-tree. The \emph{pathwidth} of $G$ is the minimum $k$ for which $G$ is a subgraph of an interval\footnote{A graph $G$ is an \emph{interval graph} if each vertex $v\in V(G)$ can be assigned an interval $I_v\subset \Re$ such that $I_w\cap I_v\neq\emptyset$ if and only if $vw\in E(V)$. } graph with no clique of order $k+2$. Note that an interval graph with no $(k+2)$-clique is a special case of a $k$-tree, and thus the treewidth of a graph is at most its pathwidth.

\lemref{trees} shows that ($1$-)trees have bounded distance-number. However, this is not true for $2$-trees since $K_{2,n}$ has treewidth (and pathwidth) at most $2$. By \thmref{degree-8}, there are $n$-vertex graphs of bounded degree with distance-number approaching $\Omega(n^{0.864138})$. On the other hand, no polynomial lower bound holds for graphs of bounded degree \emph{and} bounded treewidth, as shown in the following theorem.

%%%%%%%%%%%%%%%%%%%%%%%%
\begin{thm}\thmlabel{tw}
Let $G$ be a graph with $n$ vertices, maximum degree $\Delta$, and treewidth $k$. Then the distance-number of $G$ satisfies $$\dn{G}\in \mathcal{O}(\Delta^4 k^3\log n)\enspace.$$
\end{thm}
%%%%%%%%%%%%%%%%%%%%%%%%

To prove \thmref{tw} we use the following lemma, the proof of which is readily obtained by inspecting the proof  of Lemma~8 in \citep{DSW-CGTA}. An \emph{$H$-partition} of a graph $G$ is a partition of $V(G)$ into vertex sets $V_1, \dots, V_t$ such that $H$ is the graph with vertex set $V(H):=\{1,\dots,t\}$ where $ij\in E(H)$ if and only if there exists $v\in V_i$ and $w\in V_j$ such that $v_iv_j\in E(G)$. 
%vertices obtained from $G$ by identifying all the vertices \footnote{and deleting loops and replacing parallel edges by a single edge} in each set $V_i$, $1\leq i\leq t$. 
The \emph{width} of an $H$-partition is $\max\{|V_i|:1\leq i\leq t\}$.

%%%%%%%%%%%%%%%%%%%%%%%%
\begin{lem}[\citep{DSW-CGTA}]\lemlabel{helper}
Let $H$ be a graph admitting a drawing $D$ with $s$ distinct edge-slopes and $\ell$ distinct edge-lengths. Let $G$ be a graph admitting an $H$-partition of width $w$. Then the distance-number of $G$ satisfies 
$$\dn{G} \leq s\ell w(w-1)+\floor{\frac{w}{2}} + \ell\enspace.$$
\end{lem}
%%%%%%%%%%%%%%%%%%%%%%%%

\begin{proof}[Sketch of Proof]
The general approach is to start with $D$ and then replace each vertex of $H$ by a sufficiently scaled down and appropriate rotated copy of the drawing of $K_w$ on a regular $w$-gon. The only difficulty is choosing the rotation and the amount by which to scale the $w$-gons so that we obtain a (non-degenerate) drawing of $G$. Refer to \cite{DSW-CGTA} for the full proof.
\end{proof}

\begin{proof}[Proof of \thmref{tw}]
Let $w$ be the minimum  width of a $T$-partition of $G$ in which $T$ is a tree. The best known upper bound is $w\leq\frac{5}{2}(k+1)(\frac{7}{2}\,\Delta(G)-1)$, which was obtained by \citet{Wood-TreePartitions} using a minor improvement to a similar result by an anonymous referee of the paper by \citet{DO-JGT95}. For each vertex $x\in V(T)$, there are at most $w\Delta$ edges of $G$ incident to vertices mapped to $x$. Hence we can assume that $T$ is a forest with maximum degree $w\Delta$, as otherwise there is an edge of $T$ with no edge of $G$ mapped to it, in which case the edge of $T$ can be deleted. Similarly, $T$ has at most $n$ vertices. \citet{Scheffler89} proved that $T$ has pathwidth at most $\log(2n + 1)$; see \cite{Bodlaender-TCS98}.  \citet{DSW-CGTA} proved that every tree $T$ with pathwidth $p\geq 1$ has a drawing with $\max\{\Delta(T)-1, 1\}$ slopes and $2p- 1$ edge-lengths. Thus $T$ has a drawing with at most $\Delta w - 1$ slopes and at most $2 \log(2n + 1) - 1$ edge-lengths. By \lemref{helper}, $$\dn{G}\leq (\Delta w -1)(2\log (2n+1)-1)w(w-1) + \floor{\frac{w}{2}} + 2\log(2n + 1) - 1,$$
which is in $\mathcal{O}(\Delta w^3 \log n)\subseteq \mathcal{O}(\Delta^4 k^3 \log n)$.
\end{proof}

%%%%%%%%%%%%%%%%%%%%%%%%
\begin{cor}
Every $n$-vertex graph with bounded degree and bounded treewidth has distance-number $\mathcal{O}(\log n)$.
\end{cor}
%%%%%%%%%%%%%%%%%%%%%%%%

Since a path has a drawing with one slope and one edge-length, \lemref{helper} with $s=\ell=1$ implies that every graph $G$ with a $P$-partition of width $k$ for some path $P$ has distance-number $\dn{G}\leq k(k-\tfrac{1}{2}) + 1$.

%%%%%%%%%%%%%%%%%%%%%%%%%%%%%%%%%%%%%%%%%%%%%%%%%%%%%%%%%%%%%%%%%%%%%%%%%%%%%
\section{Bandwidth}\seclabel{Bandwidth}
%%%%%%%%%%%%%%%%%%%%%%%%%%%%%%%%%%%%%%%%%%%%%%%%%%%%%%%%%%%%%%%%%%%%%%%%%%%%%

This section establishes an upper bound on the distance-number in terms of the bandwidth. Let $G$ be a graph. A \emph{vertex ordering} of $G$ is a bijection $\sigma:V(G)\rightarrow\{1,2,\dots,|V(G)|\}$. The \emph{width} of $\sigma$ is defined to be $\max\{|\sigma(v)-\sigma(w)|:vw\in E(G)\}$. The \emph{bandwidth} of $G$, denoted by \bw{G}, is the minimum width of a vertex ordering of $G$. The \emph{cyclic width} of $\sigma$ is defined to be $\max\{\min\{|\sigma(v)-\sigma(w)|,n-|\sigma(v)-\sigma(w)|\}:vw\in E(G)\}$. The \emph{cyclic bandwidth} of $G$, denoted by \cbw{G}, is the minimum cyclic width of a vertex ordering of $G$; see \citep{CLS-DAM08, LSC-DM02, Zhou00, HKR-BICA99, Lin-SSMS94}. Clearly $\cbw{G}\leq\bw{G}$ for every graph $G$. 

%CY-JCO07, Hao01, LSC-AC97, Lin-Networks97, 

\begin{lem}
\lemlabel{bandwidth}
For every graph $G$,
$$\dn{G}\leq\cbw{G}\leq\bw{G}\enspace.$$
\end{lem}

\begin{proof}
Given a vertex ordering $\sigma$ of an $n$-vertex $G$, position the vertices of $G$ on a regular $n$-gon in the order $\sigma$. We obtain a drawing of $G$ in which the length of each edge $vw$ is determined by 
$$\min\{|\sigma(v)-\sigma(w)|,n-|\sigma(v)-\sigma(w)|\}\enspace.$$  
Thus the number of edge-lengths equals 
$$|\{\min\{|\sigma(v)-\sigma(w)|,n-|\sigma(v)-\sigma(w)|\}:vw\in E(G)\}|,$$ which is at most the cyclic width of $\sigma$. The result follows. 
\end{proof}

\begin{cor}
The distance-number of every $n$-vertex degree-$\Delta$ planar graph $G$ satisfies
$$\dn{G}\leq\frac{15n}{\log_{\Delta}n}\enspace.$$
\end{cor}

\begin{proof}
\citet{BPTW-Bandwidth-TGGT} proved that $\bw{G}\leq\frac{15n}{\log_{\Delta}n}$. The result follows from \lemref{bandwidth}.
\end{proof}

%%%%%%%%%%%%%%%%%%%%%%%%%%%%%%%%%%%%%%%%%%%%%%%%%%%%%%%%%%%%%%%%%%%%%%%%%%%%%
\section{Cartesian Products}\seclabel{CartesianProducts}
%%%%%%%%%%%%%%%%%%%%%%%%%%%%%%%%%%%%%%%%%%%%%%%%%%%%%%%%%%%%%%%%%%%%%%%%%%%%%

This section discusses the distance-number of cartesian products of graphs. For graphs $G$ and $H$, the \emph{cartesian product} $G \square H$ is the graph with vertex set $V(G\square H) := V(G) \times V(H)$,
where $(v,w)$ is adjacent to $(p,q)$ if and only if
(1) $v=p$ and $wq$ is an edge of $H$, or
(2) $w=q$ and $vp$ is an edge of $G$. 

Thus $G \square H$ is the grid-like graph with a copy of $G$ in each row and a copy of $H$ in each column. Type (1) edges form copies of $H$, and type (2) edges form copies of $G$. For example, $P_n \square P_n$ is the planar grid, and $C_n \square C_n$ is the toroidal grid. 

The cartesian product is associative and thus multi-dimensional products are well defined. For example, the $d$-dimensional product $K_2\square K_2\square \dots\square K_2$ is the $d$-dimensional hypercube $Q_d$. It is well known that $Q_d$ is a unit-distance graph. \citet{HP-Bled} proved that the cartesian product operation preserves unit-distance graphs. That is, if $G$ and $H$ are unit-distance graphs, then so is $G\square H$, as illustrated in \figref{333}. The following theorem generalises this result.

\begin{figure}[!ht]
\begin{center}
\includegraphics{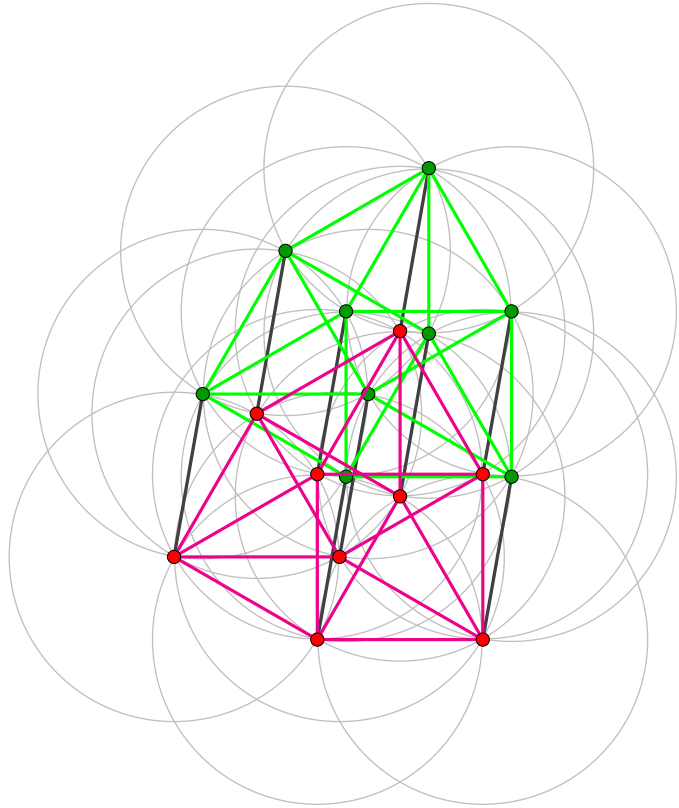}
\end{center}
\caption{A unit-distance drawing of $K_3\square K_3\square K_2$}
\figlabel{333}
\end{figure}

\begin{thm}\thmlabel{CartProd}
For all graphs $G$ and $H$, the distance-numbers of $G\square H$ satisfy
\begin{align*}
\max\{\ddn{G},\ddn{H}\} &\leq \ddn{G\square H} \leq \ddn{G}+\ddn{H}-1\enspace,\text{ and}\\
\max\{\dn{G},\dn{H}\} &\leq \dn{G\square H}\leq\dn{G}+\dn{H}-1\enspace.
\end{align*}
\end{thm}

\begin{proof}
The lower bounds follow since $G$ and $H$ are subgraphs of $G\square H$. We prove the upper bound for \dn{G\square H}. The proof for \ddn{G\square H} is simpler. 

Fix a drawing of $G$ with $\dn{G}$ edge-lengths. 
Let $(\x(v),\y(v))$ be the coordinates of each vertex $v$ of $G$ in this drawing.
Fix a drawing of $H$ with $\dn{H}$ edge-lengths, scaled so that one edge-length in the drawing of $G$ coincides with one edge-length in the drawing of $H$.
Let $\alpha$ be a real number in $[0,2\pi)$.
Let $(\x_\alpha(w),\y_\alpha(w))$ be the coordinates of each vertex $w$ of $G$ in this drawing of $H$ rotated by $\alpha$ degrees about the origin.

Position vertex $(v,w)$ in $G\square H$ at $(\x(v)+\x_\alpha(w),\y(v)+\y_\alpha(w))$. This mapping preserves edge-lengths. In particular, the length of a type-(1) edge $(v,u)(v,w)$ equals the length of the edge $uw$ in $H$, and the length of a type-(2) edge $(u,v)(w,v)$ equals the length of the edge $uw$ in $G$. Thus for each $\alpha$, the mapping of $G\square H$ has $\dn{G}+\dn{H}-1$ edge-lengths. 

It remains to prove that for some $\alpha$ the mapping of $G\square H$ is a drawing. That is, no vertex intersects the closure of an incident edge. An angle $\alpha$ is \emph{bad} for a particular vertex/edge pair of $G\square H$ if that vertex intersects the closure of that edge in the mapping with rotation $\alpha$. 

Observe that the trajectory of a vertex $(v,w)$ of $G\square H$ (taken over all $\alpha$) is a circle centred at $(\x(v),\y(v))$ with radius $\dist_H(0,w)$. 

Now for distinct points $p$ and $q$ and a line $\ell$, there are only two angles $\alpha$ such that the rotation of $p$ around $q$ by an angle of $\alpha$ contains $\ell$ (since the trajectory of $p$ is a circle that only intersects $\ell$ in two places), and there are only two angles $\alpha$ such that the rotation of $\ell$ around $q$ by an angle of $\alpha$ contains $p$.

It follows that there are finitely many bad values of $\alpha$ for a particular vertex/edge pair of $G\square H$. Hence  there are finitely many bad values of $\alpha$ in total. Hence some value of $\alpha$ is not bad for every vertex/edge pair in $G\square H$. Hence $D_\alpha$ is a valid drawing of $G\square H$. 
\end{proof}

Note that \citet{LT-NM66} proved a result analogous to \thmref{CartProd} for dimension. 

Let $G^d$ be the $d$-fold cartesian product of a graph $G$. The same construction used in \thmref{CartProd} proves the following:

\begin{thm}\thmlabel{CartProdSimple}
For every graph $G$ and integer $d\geq1$, the distance-numbers of $G^d$ satisfy
$$\ddn{G^d}=\ddn{G}\quad\text{ and }\quad\dn{G^d}=\dn{G}\enspace.$$
\end{thm}

%\begin{figure}[H]\begin{center}\includegraphics{33}\end{center}\caption{A unit-distance drawing of $K_3\square K_3$}\end{figure}

%\begin{figure}[H]\begin{center}\includegraphics{333}\end{center}\caption{Oh, what a tangled web we weave: a unit-distance drawing of $K_3\square K_3\square K_3$, I think!}\end{figure}

%%%%%%%%%%%%%%%%%%%%%%%%%%%%%%%%%%%%%%%%%%%%%%%%%%%%%%%%%%%%%%%%%%%%%%%%%%%%%
\section{Open Problems}\seclabel{Conclusion}
%%%%%%%%%%%%%%%%%%%%%%%%%%%%%%%%%%%%%%%%%%%%%%%%%%%%%%%%%%%%%%%%%%%%

We conclude by mentioning some of the many open problems related to distance-number. 

\begin{itemize}

\item What is $\dn{K_n}$? We conjecture that $\dn{K_n}=\floor{\frac{n}{2}}$. That is, every set of $n$ points in general position determine at least $\floor{\frac{n}{2}}$ distinct distances. Note that \citet{Altman-AMM63,Altman-CMB72} proved this conjecture for points in convex position.

\item What is the relationship between distance-number and degenerate distance-number? In particular, is there a function $f$ such that $\dn{G}\leq f(\ddn{G})$ for every graph $G$? 

\item \threethmref{degree-7}{degree-8}{degree-5} establish a lower bound for the distance-number of bounded degree graphs. But no non-trivial upper bound is known. Do $n$-vertex graphs with bounded degree have distance-number in $o(n)$? 

\item Outerplanar graphs have distance-number in $\mathcal{O}(\Delta^4\log n)$ by \thmref{tw}. Do outerplanar graphs (with bounded degree) have bounded (degenerate) distance-number? 

\item Non-trivial lower and upper bounds on the distance-numbers are not known for many other interesting graph families including: degree-$3$ graphs, degree-$4$ graphs, $2$-degenerate graphs with bounded degree, graphs with bounded degree and bounded pathwidth.  

\item As described in \secref{Motivation}, determining the maximum number of times the unit-distance can appear among $n$ points in the plane is a famous open problem. We are unaware if the following apparently simpler tasks have been attempted: Determine the maximum number of times the unit-distance can occur among $n$ points in the plane such that no three are collinear. Similarly, determine the maximum number of edges in an $n$-vertex graph $G$ with $\dn{G}=1$.

\item Determining the maximum chromatic number of unit-distance graphs in $\Re^d$ is a well-known open problem. The best known upper bound of $(3+o(1))^d$ is due to \citet{LR-Math72}. Exponential lower bounds are known \citep{Raig-UMN00,FW-Comb81}. Unit-distance graphs in the plane are $7$-colourable \citep{HD64}, and thus $\chi(G)\leq 7^{\ddn{G}}$. Can this bound be improved? 

\item Degenerate distance-number is not bounded by any function of dimension since $K_{n,n}$ has dimension $4$ and unbounded degenerate distance-number. On the other hand, $\DIM(G)\leq 2\cdot\chi(G)\leq 2\cdot 7^{\ddn{G}}$. Is $\DIM(G)$ bounded by a polynomial function of $\ddn{G}$? 

\item Every planar graph has a crossing-free drawing. A long standing open problem involving edge-lengths, due to Harborth~\etal~\citep{MR1830610, MR1126238, MR905393}, asks whether every planar graph has a crossing-free drawing in which the length of every edge is an integer. \citet{GGM-JGT08} recently answered this question in the affirmative for cubic planar graphs. \citet{Archdeacon-Integer} extended this question to nonplanar graphs and asked what is the minimum $d$ such that a given graph has a crossing-free drawing in $\Re^d$ with integer edge-lengths. Note that every $n$-vertex graph has such a drawing in $\Re^{n-1}$.

\item The \emph{slope number} of a graph $G$, denoted by \sn{G}, is the minimum number of edge-slopes over all drawings of $G$. \citet{DESW-CGTA} established results concerning the slope-number of planar graphs. \citet{KPPT-GD07} proved that degree-$3$ graphs have slope-number at most $5$. On the other hand, \citet{BMW-EJC06} and \citet{PachPal-EJC06} independently proved that there are $5$-regular graphs with arbitrarily large slope number. Moreover, for $\Delta\geq7$, \citet{DSW-CGTA} proved that there are $n$-vertex degree-$\Delta$ graphs whose slope number is at least $n^{1-\frac{\varepsilon}{\Delta+4}}$. The proofs of these results are similar to the proofs of \threethmref{degree-7}{degree-8}{degree-5}. Given that \thmref{tw} also depends on slopes, it it tempting to wonder if there is a deeper connection between slope-number and distance-number. For example, is there a function $f$ such that $\sn{G} \leq f(\Delta(G),\dn{G})$ and/or $\dn{G} \leq f(\sn{G})$ for every graph $G$. Note that some dependence on $\Delta(G)$ is necessary since $\sn{K_{1,n}}\rightarrow\infty$ but $\dn{K_{1,n}}=1$.

\end{itemize}

%\bibliographystyle{myNatbibStyle}
%\bibliography{myBibliography,myConferences}

\def\soft#1{\leavevmode\setbox0=\hbox{h}\dimen7=\ht0\advance \dimen7
  by-1ex\relax\if t#1\relax\rlap{\raise.6\dimen7
  \hbox{\kern.3ex\char'47}}#1\relax\else\if T#1\relax
  \rlap{\raise.5\dimen7\hbox{\kern1.3ex\char'47}}#1\relax \else\if
  d#1\relax\rlap{\raise.5\dimen7\hbox{\kern.9ex \char'47}}#1\relax\else\if
  D#1\relax\rlap{\raise.5\dimen7 \hbox{\kern1.4ex\char'47}}#1\relax\else\if
  l#1\relax \rlap{\raise.5\dimen7\hbox{\kern.4ex\char'47}}#1\relax \else\if
  L#1\relax\rlap{\raise.5\dimen7\hbox{\kern.7ex
  \char'47}}#1\relax\else\message{accent \string\soft \space #1 not
  defined!}#1\relax\fi\fi\fi\fi\fi\fi} \def\Dbar{\leavevmode\lower.6ex\hbox to
  0pt{\hskip-.23ex \accent"16\hss}D}

\end{document}